\documentclass[12pt,leqno,a4paper]{annacad}
\pdfoutput=1
\usepackage[T1]{fontenc}
\usepackage[utf8]{inputenc}
\usepackage[british,finnish]{babel}
\usepackage{lmodern}
\usepackage{microtype}
\usepackage{mathtools}
\usepackage{esint}
\usepackage{comment}
\usepackage{hyperref}
\usepackage{amssymb}
\usepackage{textcomp}
\usepackage[toc]{appendix}
\usepackage{pdfpages}

\title{Calder\'on's problem for $p$-Laplace type equations}
\author{Tommi Brander}

\newcommand*{\R}{\mathbb{R}}

\newcommand*{\N}{\mathbb{N}}
\newcommand*{\Z}{\mathbb{Z}}

\newcommand*{\eps}{\varepsilon}
\newcommand*{\doo}{\partial}
\newcommand*{\ol}[1]{\overline{#1}}
\newcommand{\sulut}[1]{\left( #1 \right)}
\newcommand{\joukko}[1]{\left\{ #1 \right\}}
\newcommand{\abs}[1]{\left\lvert #1 \right\rvert}

\newcommand{\der}{\mathrm{d}}

\newcommand{\ip}[2]{\left\langle#1,#2\right\rangle}

\theoremstyle{plain}
\newtheorem{theorem}{Theorem}
\newtheorem{lemma}[theorem]{Lemma}

\theoremstyle{definition}
\newtheorem{definition}[theorem]{Definition}

\theoremstyle{remark}
\newtheorem{remark}[theorem]{Remark}

\DeclareMathOperator{\dive}{div}

\newcommand{\enquote}[1]{`#1'}

\addto\captionsbritish{%
  \renewcommand{\contentsname}{Contents of the introductory part}
}

\def\vint_#1{\mathchoice%
          {\mathop{\kern 0.2em\vrule width 0.6em height 0.69678ex
depth -0.58065ex
                  \kern -0.8em \intop}\nolimits_{\kern -0.4em#1}}%
          {\mathop{\kern 0.1em\vrule width 0.5em height 0.69678ex
depth -0.60387ex
                  \kern -0.6em \intop}\nolimits_{#1}}%
          {\mathop{\kern 0.1em\vrule width 0.5em height 0.69678ex
              depth -0.60387ex
                  \kern -0.6em \intop}\nolimits_{#1}}%
          {\mathop{\kern 0.1em\vrule width 0.5em height 0.69678ex
depth -0.60387ex
                  \kern -0.6em \intop}\nolimits_{#1}}}
\def\vintslides_#1{\mathchoice%
          {\mathop{\kern 0.1em\vrule width 0.5em height 0.697ex depth -0.581ex
                  \kern -0.6em \intop}\nolimits_{\kern -0.4em#1}}%
          {\mathop{\kern 0.1em\vrule width 0.3em height 0.697ex depth -0.604ex
                  \kern -0.4em \intop}\nolimits_{#1}}%
          {\mathop{\kern 0.1em\vrule width 0.3em height 0.697ex depth -0.604ex
                  \kern -0.4em \intop}\nolimits_{#1}}%
          {\mathop{\kern 0.1em\vrule width 0.3em height 0.697ex depth -0.604ex
                  \kern -0.4em \intop}\nolimits_{#1}}}

\newcommand{\aveint}[2]{\mathchoice%
          {\mathop{\kern 0.2em\vrule width 0.6em height 0.69678ex
depth -0.58065ex
                  \kern -0.8em \intop}\nolimits_{\kern -0.45em#1}^{#2}}%
          {\mathop{\kern 0.1em\vrule width 0.5em height 0.69678ex
depth -0.60387ex
                  \kern -0.6em \intop}\nolimits_{#1}^{#2}}%
          {\mathop{\kern 0.1em\vrule width 0.5em height 0.69678ex
depth -0.60387ex
                  \kern -0.6em \intop}\nolimits_{#1}^{#2}}%
          {\mathop{\kern 0.1em\vrule width 0.5em height 0.69678ex
depth -0.60387ex
                  \kern -0.6em \intop}\nolimits_{#1}^{#2}}}

\numberwithin{theorem}{section}
\numberwithin{equation}{section}

\begin{document}
\pagenumbering{gobble}
\selectlanguage{british}


\includepdf[pages={3-4}]{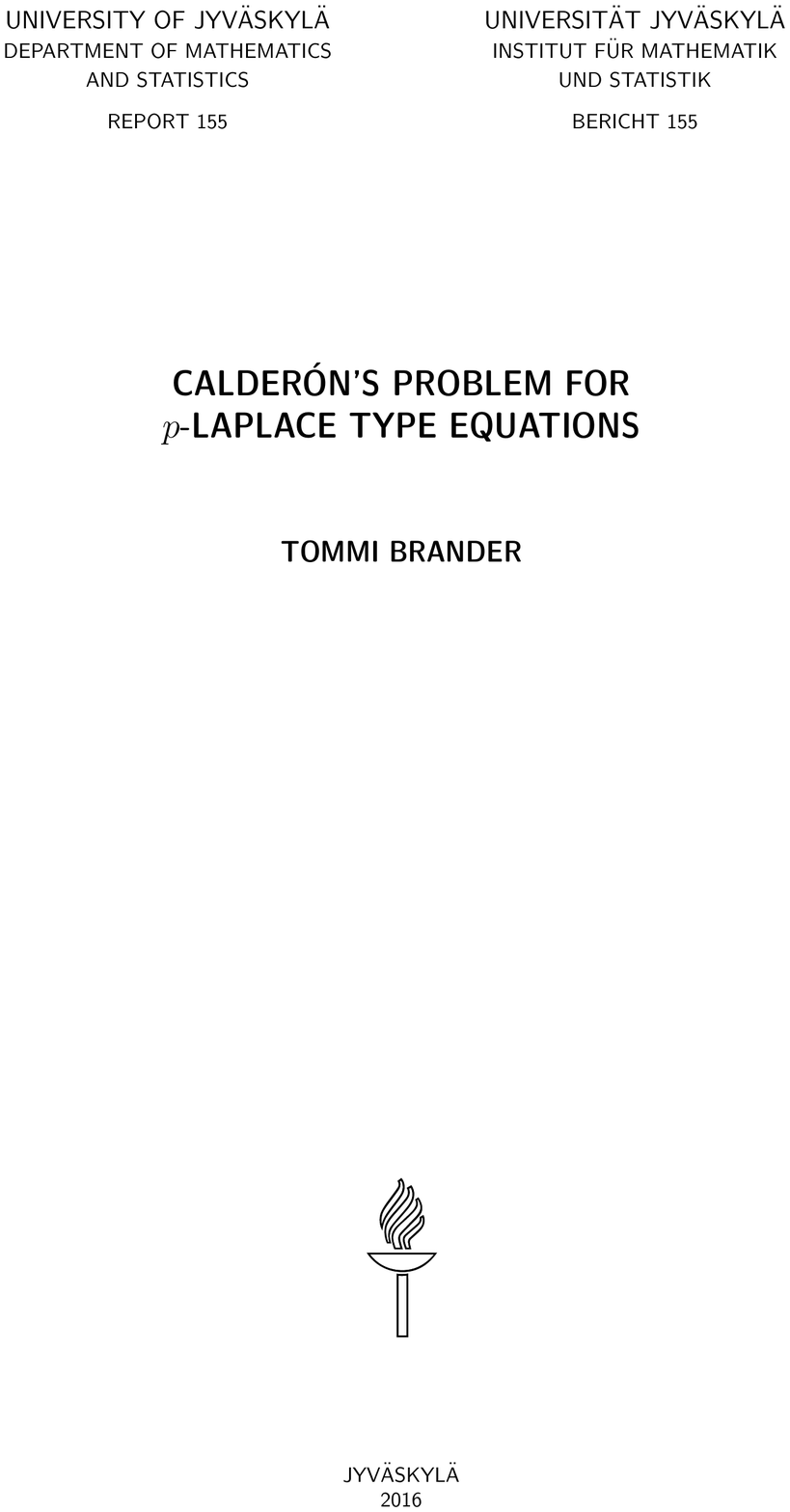} %

\newpage

\pagenumbering{roman}

\section*[Abstract]{Abstract}

We investigate a generalization of Calder\'on's problem of recovering the conductivity coefficient in a conductivity equation from boundary measurements.
As a model equation we consider the $p$-conductivity equation
\begin{equation*}
\dive \sulut{\sigma \abs{\nabla u}^{p-2} \nabla u} = 0
\end{equation*}
with $1 < p < \infty$, which reduces to the standard conductivity equation when $p = 2$.

The thesis consists of results on the direct problem, boundary determination and detecting inclusions.
We formulate the equation as a variational problem also when the conductivity~$\sigma$ may be zero or infinity in large sets.
As a  boundary determination result we recover the first order derivative of a smooth conductivity on the boundary.
We use the enclosure method of Ikehata to recover the convex hull of an inclusion of finite conductivity and find an upper bound for the convex hull if the conductivity within an inclusion is zero or infinite.

\newpage

\begin{otherlanguage}{finnish}

\section*[Tiivistelmä]{Tiivistelmä}

Calder\'onin ongelma kysyy: voidaanko johtavuuden arvo johtavuusyhtälössä määrittää reunamittauksia käyttäen?
Tutkimme Calder\'onin ongelman yleistystä tilanteeseen, jossa suoraa ongelmaa kuvaa $p$-johtavuusyhtälö
\begin{equation*}
\dive \sulut{\sigma \abs{\nabla u}^{p-2} \nabla u} = 0,
\end{equation*}
missä $1 < p < \infty$.
Yhtälö on tavallinen johtavuusyhtälö, jos $p = 2$.

Väitöskirjan tulokset koskevat suoraa ongelmaa, reunamääritystä ja sisältymän (eli inkluusion tai esteen) etsimistä.
Yleistämme suoran ongelman tilanteeseen, jossa johtavuus~$\sigma$ voi saada arvon nolla tai ääretön suurissa joukoissa.
Reunamääritystuloksena määritämme siistin johtavuuden gradientin tutkitun alueen reunalla.
Käytämme Ikehatan kotelointimenetelmää määrittääksemme äärellisjohtavuuksisen sisältymän kuperan verhon, ja ylärajan verholle, jos johtavuus sisältymässä on nolla tai ääretön.

\end{otherlanguage}

\selectlanguage{british}

\newpage

\section*[Acknowledgements]{Acknowledgements}
I dedicate this thesis to the memory of Juha Brander.

I would like to to thank my advisor, professor Mikko Salo, for his guidance and help, professor Sergey Repin for showing me that research is a possible path for me, and fellow students and faculty for an encouraging environment.

The research was partially funded by Academy of Finland through the Finnish Centre of Excellence in Inverse Problems Research, for which I am grateful.

Finally, most of my gratitude belongs to Terhi for everything and to Noora for all the joy she brought to us.

\newpage

\section*[List of included articles]{List of included articles}

\noindent
This thesis is based on the work contained within the following publications:

\vspace{1em}
\begingroup
\renewcommand{\section}[2]{}

\endgroup

\noindent
The author has participated actively in development of the joint papers~\cite{Brander:Kar:Salo:2015:thesis,Brander:Ilmavirta:Kar:2015:thesis}.

\newpage

\tableofcontents

\newpage\null\thispagestyle{empty}\newpage

\pagenumbering{arabic}
\section{Inverse problems}

The field of inverse problems is motivated by practical problems, where the objective is to recover information on a medium by indirect measurements; for example:
\begin{itemize}
\item To what extent do X-rays penetrate different human tissue and how much of them is absorbed?
	After measuring the absorption of X-rays, the inverse problem is to calculate which materials the rays passed through.
\item How does the composition of the crust of Earth affect the movement of shock waves and sound waves?
	The inverse problem is to use measurements of the waves to find out the material properties of Earth.
\item How does the value of an option (a financial instrument) depend on its volatility?
	The inverse problem is to determine the volatility of an option from the prices with which it is bought and sold in the free market.
\item How does the heat conductivity of an object affect the distribution of heat within it?
The inverse problem is to determine the heat conductivity by measuring heat and heat flow on the surface of the object.
\item Electrical impedance tomography asks the previous question, but with heat replaced by electricity.
\end{itemize}
The inverse problems examined in this thesis are related to recovering the conductivity of heat or electricity from surface measurements.

\section{Calder\'on's problem and electrical impedance tomography}

We investigate $p$-Calder\'on's problem in one dimension before posing the standard Calder\'on's problem in higher dimensions; in our nomenclature, the standard problem is called $2$-Calder\'on's problem.

We use the following notations throughout the thesis:
\begin{itemize}
\item Open bounded set~$\Omega \subset \R^d$, where $d \in \Z_+$ is the dimension.
We assume~$d \geq 2$ except in section~\ref{sec:1d}, where~$d=1$.
\item Potential~$u \colon \Omega \to \R$.
The potential solves the $p$-conductivity equation, which will be discussed below.
The potential may be electric potential (electric field potential or electrostatic potential), or temperature.
We only consider steady state equations, where the potential does not depend on time.
\item Conductivity~$\sigma \colon \Omega \to [0,\infty]$, which we always assume to be measurable.
The conductivity may be electrical conductivity or heat conductivity.
Sometimes we consider an extension to $\ol \Omega$, which we also write as $\sigma$.
\end{itemize}
Also, in the following physical introduction only, we use the current flux or heat flux~$j \colon \Omega \to \R^d$.

\subsection{Physical backgrond} \label{subsec:physics}

The negative current flux is proportional to differences in potential and to conductivity by Ohm's law or by Fourier's law of thermal conduction:
\begin{equation} \label{eq:fys1}
-j = \sigma \nabla u
\end{equation}
We assume there are no sources or sinks of electricity or heat inside the object~$\Omega$, whence by conservation of energy (or Gauss's law or Kirchhoff's law)
\begin{equation} \label{eq:fys2}
\dive j = 0.
\end{equation}

Thus, by equations~\eqref{eq:fys1} and \eqref{eq:fys2}, the conductivity equation
\begin{equation}
-\dive \sulut{\sigma \nabla u} = 0
\end{equation}
models the steady state of conduction of heat or electricity.

The thesis discusses a non-linear variant of the conductivity equation, where the linear relation~\eqref{eq:fys1} is replaced with a non-linear power law relation 
\begin{equation} \label{eq:fys3}
-j = \sigma \abs{\nabla u}^{p-2} \nabla u,
\end{equation}
which leads to the $p$-conductivity equation
\begin{equation}
-\dive \sulut{\sigma \abs{\nabla u}^{p-2} \nabla u} = 0.
\end{equation}
The materials that do not obey Ohm's law are called non-Ohmic.
The power-law behaviour is one possible non-Ohmic behaviour and may be a fine approximation in more complicated non-Ohmic situations.
In context of electricity, power-law behaviour has been observed in certain polycrystalline materials near the superconducting-normal transition~\cite{Dubson:Herbert:Calabrese:Harris:Patton:Garland:1988,Bueno:Longo:Varela:2008}.
Other applications of the $p$-Laplace equation, which is the $p$-conductivity equation with~$\sigma \equiv 1$, include nonlinear dielectrics~\cite{Bueno:Longo:Varela:2008,Garroni:Kohn:2003,Garroni:Nesi:Ponsiglione:2001,Kohn:Levy:1998,Talbot:Willis:1994:a,Talbot:Willis:1994:b} and plastic moulding~\cite{Aronsson:1996}. It is also used to model electro-rheological and thermo-rheological fluids~\cite{Antontsev:Rodrigues:2006,Berselli:Diening:Ruzicka:2010,Ruzicka:2000}, fluids governed by a power law~\cite{Aronsson:Janfalk:1992}, viscous flows in glaciology~\cite{Glowinski:Rappaz:2003} and some plasticity phenomena~\cite{Atkinson:Champion:1984,Idiart:2008,PonteCastaneda:Suquet:1998,PonteCastaneda:Willis:1985,Suquet:1993}.
There are further applications to image processing~\cite{Kuijper:2007} and conformal geometry~\cite{Liimatainen:Salo:2012,Julin:Liimatainen:Salo:2015}.

\subsection{One dimensional case}\label{sec:1d}

As a simple special case we investigate a one-dimensional situation, where the domain $\Omega$ is an interval, which we write as the open\footnote{I use the customary Finnish notation for open interval, $]a,b[$, rather than the notation $\sulut{a,b}$, since the first is clearly distinct from a vector in dimension~two.} interval $]a,b[$ with $a,b \in \R$.
The treatment is somewhat heuristical, as the precise function space where the problem is solved is not defined here, the variational formulation for the forward problem is not justified and the definition of the Dirichlet to Neumann map is not justified.
We also use the strong formulation of the forward problem without any worries concerning regularity when doing so is expedient.
We give a rigorous treatment in article~\cite{Brander:Ilmavirta:Kar:2015:thesis}; see also section~\ref{sec:direct}.

Here and elsewhere the parameter~$p$ may always take values strictly between one and infinity, unless otherwise mentioned.
The one-dimensional case with $p=2$ can be found at least in the unpublished manuscript~\cite{Feldman:Salo:Uhlmann:2012}.

We consider the conductivity to be given by a measurable function $\sigma \colon ]a,b[ \to [0,\infty]$, which we assume to be essentially bounded from above and away from zero in the complement of $D_0\cup D_\infty$, where and henceforth we write $D_j = \sigma^{-1}(\joukko{j})$.
We also assume that $D_0$ and $D_\infty$ are unions of finite numbers of open intervals, the closures of which are disjoint from each other and the boundary $\{ a,b \}$.

Given some Dirichlet boundary values $\sulut{A,B} \in \R^2$, the Dirichlet problem for the one-dimensional $p$-conductivity equation is
\begin{equation}  \label{eq:1-d}
\begin{cases}
\sulut{\sigma(x)\abs{u'(x)}^{p-2}u'(x)}' = 0 \text{ for } a < x < b\\
u(a) = A \\
u(b) = B.
\end{cases}
\end{equation}
The set of both Dirichlet and Neumann boundary values is $\R^2$, so the strong Dirichlet to Neumann map is $\Lambda_\sigma \colon \R^2 \to\R^2$,
\begin{equation}
\sulut{A,B} \mapsto \sulut{-\sigma(a)\abs{u'(a)}^{p-2}u'(a), \sigma(b)\abs{u'(b)}^{p-2}u'(b)},
\end{equation}
where $u$ solves the $p$-conductivity equation~\eqref{eq:1-d} and when $u$ and $\sigma$ are smooth enough.
The weak definition is $\Lambda_\sigma  \colon \R^2 \to \joukko{L \colon \R^2 \to \R; L \text{ linear}}$,
\begin{equation}
\begin{split}
&\ip{\Lambda_\sigma (\sulut{A,B})}{\sulut{\alpha,\beta}} = \int_a^b \sigma \abs{u'}^{p-2} u' h'  \der x,
\end{split}
\end{equation}
where $h \in W^{1,p}(\Omega)$, $h(a) = \alpha$, $h(b) = \beta$, $h'|_{D_\infty} = 0$, and $u$ solves the problem~\eqref{eq:1-d}.
Supposing the conductivity is somewhat regular near $a$ and $b$, the strong Dirichlet to Neumann map can be recovered from the weak one by using test functions~$h$ with
\begin{equation}
\begin{cases}
h'(x) = 1/\eps &\text{ while } x \geq b-\eps \\
h'(x) = 0 &\text{ otherwise}
\end{cases}
\end{equation}
and similar test functions with slope $-1/\eps$ near $a$ that are constant elsewhere.

Before proceeding further with the inverse problem, we consider the issues of existence and uniqueness for the forward problem~\eqref{eq:1-d}.
We find out that the correct space for the solutions is a close relative of the space $W^{1,p}(\Omega)$ with the correct Dirichlet boundary values.
For rigorous and more precise treatment of the issues see article~\cite{Brander:Ilmavirta:Kar:2015:thesis} and section~\ref{sec:direct}.

The solution is discovered by minimizing the energy
\begin{equation} \label{eq:1-d_energy}
I(v) = \int_a^b \sigma(x) \abs{v'(x)}^p \der x.
\end{equation}

We shall explicitly solve the problem~\eqref{eq:1-d} by using both the strong formulation of the problem and the energy minimization formulation~\eqref{eq:1-d_energy}.

We first  assume that $D_0 \neq \emptyset$.
A solution for our problem is
\begin{equation}
u(x) =
\begin{cases}
A \text{ for } x \leq \inf D_0 \\
g(x) \text{ for } \inf D_0  \leq x \leq \sup D_0\\
B \text{ for } x \geq \sup D_0,
\end{cases}
\end{equation}
where $g(\inf D_0) = A$, $g(\sup D_0) = B$, and $g' = 0$ outside $D_0$.
We may take $g$ to be smooth.
We have $I(u) = 0$, so $u$ minimizes the energy.
Any other minimizer must be constant outside $D_0$, since its energy must vanish.
In particular, $\Lambda_\sigma \sulut{\sulut{A,B}} = 0$.

Next we consider the more interesting situation $D_0 = \emptyset$.
By the fundamental theorem of calculus we calculate the solution of problem~\eqref{eq:1-d} to be
\begin{equation}
u(x) = A + \frac{B-A}{\int_a^b \sigma^{1/(1-p)}(t) \der t}\int_a^x \sigma^\frac{1}{1-p}(t) \der t.
\end{equation}
This makes sense even when $\sigma = \infty$; the integrand is then zero and thence the solution $u$ is locally constant in regions of infinite conductivity.

Now we have explicitly solved the forward problem:
\begin{theorem}
Suppose $\Omega \subset \R$ is an interval, which we write as the open interval $]a,b[$, and $a,b \in \R$.
Let $\sigma$ be a measurable function mapping $]a,b[$ to $[0,\infty]$.
We assume it is essentially bounded from above and away from zero in the complement of $D_0\cup D_\infty$.
We also assume that $D_0$ and $D_\infty$ are unions of finite numbers of open intervals and $\Omega \setminus D_\infty$ is a set of positive Lebesgue measure.

Then the function $u \in W^{1,p}(\Omega)$ defined as
\begin{equation}
u =
\begin{cases}
w_0 &\text{ if } D_0 \neq \emptyset \\
w_1 &\text { if } D_0 = \emptyset,
\end{cases}
\end{equation}
with
\begin{equation}
w_0 (x) =
\begin{cases}
A& \text{ for } x \leq \inf D_0 \\
g(x)& \text{ for } \inf D_0  \leq x \leq \sup D_0\\
B& \text{ for } x \geq \sup D_0,
\end{cases}
\end{equation}
where $g(\inf D_0) = A$, $g(\sup D_0) = B$, $g'|_{\Omega \setminus D_0} = 0$ and $g \in W^{1,p}(\Omega)$,
and
\begin{equation}
w_1 (x) =A + \frac{B-A}{\int_a^b \sigma^{1/(1-p)}(t) \der t}\int_a^x \sigma^\frac{1}{1-p}(t) \der t,
\end{equation}
solves the one-dimensional $p$-conductivity equation~\eqref{eq:1-d} in the weak sense with Dirichlet boundary values $(A,B)$.
We interpret $\infty^{1/\sulut{1-p}}$ as zero.
\end{theorem}

The strong Dirichlet to Neumann map is
\begin{equation}
\Lambda_\sigma (A,B) = \sulut{\frac{A-B}{\int_a^b \sigma^{1/(1-p)}(t) \der t},\frac{B-A}{\int_a^b \sigma^{1/(1-p)}(t) \der t}},
\end{equation}
since
\begin{equation}
\sigma \abs{u'}^{p-2}u' = \frac{B-A}{\int_a^b \sigma^{1/(1-p)}(t) \der t}.
\end{equation}
Since $A-B$ is known from the Dirichlet data, all we can learn is the quantity
\begin{equation}
\int_a^b \sigma^{1/(1-p)}(t) \der t.
\end{equation}

We can slightly weaken our assumptions by accepting that either $a \in \ol{D_\infty}$ or $b \in \ol{D_\infty}$ with no loss in results, but if both are true, then that is all we can observe.

In a similar way the weak Dirichlet to Neumann map gives
\begin{equation}
\begin{split}
\ip{\Lambda_\sigma \sulut{A,B}}{\sulut{\alpha,\beta}}
&= \int_a^b \sigma \abs{u'}^{p-2}u'h' \der x \\
&= \frac{B-A}{\int_a^b \sigma^{1/(1-p)}(t) \der t} \int_a^b h' \der x \\
&= \frac{\sulut{B-A}\sulut{\beta-\alpha}}{\int_a^b \sigma^{1/(1-p)}(t) \der t}
\end{split}
\end{equation}
and thus provides the same information as the strong map.

We have proven the following theorem:

\begin{theorem}[$p$-Calder\'on's problem in one dimension]
Suppose the domain $\Omega \subset \R$ is a bounded open interval $]a,b[$.
Suppose the conductivity is a measurable function $\sigma \colon ]a,b[ \to [0,\infty]$ that is essentially bounded from above and away from zero in the complement of $D_0\cup D_\infty$.
We also assume that $D_0$ and $D_\infty$ are unions of finite numbers of open intervals.
Further, we assume that the sets $\ol{D_0}$, $\ol{D_\infty}$, and $\doo \Omega$ are disjoint.

We can recover the following information from the Dirichlet to Neumann map:
\begin{enumerate}
\item If $D_0 \neq \emptyset$, then we know that this is the case, but can say nothing more.
\item If $D_0 = \emptyset$, then we recover the quantity
\begin{equation}
\int_a^b \sigma^{1/(1-p)} \der t.
\end{equation}
\end{enumerate}
\end{theorem}

\subsection{Electrical impedance tomography}

Alberto Pedro Calder\'on considered a method for detecting oil wells by electrical measurements in the article~\enquote{On an inverse boundary value problem}~\cite{Calderon:1980}, which initiated the modern mathematical study of electrical impedance tomography.
The paper was published in 1980 and reprinted in~2006~\cite{Calderon:2006}.

The physical principles were outlined in section~\ref{subsec:physics}.
We now consider the equation
\begin{equation}
-\dive \sulut{\sigma \nabla u} = 0
\end{equation}
in a bounded domain $\Omega \subset \R^d$, where $u$ is the electrical potential (or voltage) and $\sigma \nabla u$ is the current flux.
The Dirichlet data $u|_{\doo \Omega}$ corresponds to boundary measurements of the potential and the Neumann data $\sigma \nabla u \cdot \nu |_{\doo \Omega}$ corresponds to current flux density across the boundary.

The idea of electrical impedance tomography is to prescribe either the Dirichlet or the Neumann data and measure the other.
Since both the Dirichlet and the Neumann boundary value problems for the conductivity equation have a unique solution for any applicable boundary data\footnote{In case of Neumann data, the solution is unique up to constants.}, the Dirichlet to Neumann map $\Lambda_\sigma$ and the Neumann to Dirichlet map, which we do not use in this thesis, are well-defined.
The strong formulation for the Dirichlet to Neumann map is
\begin{equation}
\Lambda_\sigma (f) = \sigma \nabla u \cdot \nu|_{\doo \Omega},
\end{equation}
where $u$ solves the conductivity equation with Dirichlet boundary data~$f$.
For the rigorous weak formulation of the Dirichlet to Neumann map we refer to section~\ref{sec:dn}.

At least the following kinds of results are known for Calder\'on's problem:
\begin{enumerate}
\item boundary uniqueness
\item boundary reconstruction
\item higher order derivatives at the boundary
\item stability at the boundary
\item recovering inclusions; of finite positive, zero or infinite conductivity
\item detecting cracks
\item interior uniqueness
\item interior reconstruction
\item interior stability
\item numerical algorithms
\end{enumerate}
The unpublished manuscript of Feldman, Salo and Uhlmann~\cite{Feldman:Salo:Uhlmann:2012} is a good introduction.
Uhlmann~\cite{Uhlmann:2009} has written a survey of the $2$-Calder\'on's problem.
For review of more numerical nature we refer to Borcea~\cite{Borcea:2002}.

At least the following kinds of variants or generalizations have been considered:
\begin{enumerate}
\item linearized problem
\item partial data
\item anisotropic problem
\item different base equation or system of equations
\end{enumerate}

This thesis concerns the final kind of generalization; we consider the quasilinear $p$-conductivity equation or weighted $p$-Laplace equation, which we discuss in section~\ref{sec:direct}.
Most of the results mentioned above are open for $p$-conductivity equations.

\section{$p$-conductivity equation} \label{sec:direct}

The Dirichlet problem for the linear conductivity equation, $\dive \sulut{\sigma \nabla u} = 0$, is solved by minizing the energy functional
\begin{equation}
v \mapsto \int_\Omega \sigma \abs{\nabla v}^2 \der x
\end{equation}
in the Sobolev space $W_f^{1,2}\sulut{\Omega} = f + W_0^{1,2}\sulut{\Omega}$, where the Dirichlet boundary values are $f$.

The $p$-conductivity equation generalizes the linear conductivity equation much as the $p$-Laplace equation generalizes the linear Laplace equation.
Let $1 < p < \infty$.
The $p$-conductivity equation is $\dive \sulut{\sigma \abs{\nabla u}^{p-2} \nabla u} = 0$.
The Dirichlet problem is solved by minimizing the energy
\begin{equation}
v \mapsto \int_\Omega \sigma \abs{\nabla v}^p \der x
\end{equation}
in the Sobolev space $W_f^{1,p}\sulut{\Omega} = f + W_0^{1,p}\sulut{\Omega}$, where the Dirichlet boundary values are~$f$.
The quasilinear $p$-Laplace equation, which is the special case of $p$-conductivity equation where $\sigma \equiv 1$, has been widely studied.
A good introduction is provided by the lecture notes of Peter Lindqvist~\cite{Lindqvist:2006}.
The even more general $A$-harmonic functions and related quasilinear elliptic partial differential equations of second order have been widely covered in various monographs~\cite{Ladyzhenskaya:Ural'tseva:1968, Gilbarg:Trudinger:1983,Heinonen:Kilpelainen:Martio:1993}.

By the direct method of calculus of variations, the equation is well-defined for measurable $\sigma \colon \Omega \to [0,\infty]$ bounded from above and below by positive constants.
In article~\cite{Brander:Ilmavirta:Kar:2015:thesis} the conductivity can take values zero and infinity in large sets:
\begin{theorem} \label{thm:direct}
Let $\Omega\subset\R^n$ be a bounded domain and let $\sigma\colon\Omega\to[0,\infty]$ be a measurable function.
We write $D_0=\sigma^{-1}(\joukko{0})$ and $D_\infty=\sigma^{-1}(\joukko{\infty})$ and suppose the following:
\begin{itemize}
\item The sets~$D_0$ and~$D_\infty$ are open.
\item The sets~$\ol D_0$, $\ol D_\infty$ and~$\partial\Omega$ are disjoint.
\item The conductivity~$\sigma$ is essentially bounded from below and above by positive constants in $\Omega\setminus (D_0\cup D_\infty)$.
\item The set~$D_0$ has Lipschitz boundary.
\end{itemize}
Fix $p\in]1,\infty[$.
Given any $f\in W^{1,p}(\Omega)$, there is a minimizer $u\in f+W^{1,p}_0(\Omega)$ to the energy
\begin{equation}
E(u)
=
\int_\Omega\sigma\abs{\nabla u}^p.
\end{equation}
The minimal energy is finite and the minimizer is unique up to functions that have zero gradient outside~$D_0$ and zero Dirichlet boundary values on~$\doo \Omega$.
The minimizer satisfies $\dive(\sigma\abs{\nabla u}^{p-2}\nabla u)=0$ in $\Omega\setminus(\ol D_0\cup\ol D_\infty)$ and $\nabla u=0$ in~$D_\infty$ in the weak sense.
\end{theorem}

The solutions are locally constant in $D_\infty$ and can take essentially arbitrary values in $D_0$.
Elsewhere they satisfy the usual $p$-conductivity equation in the weak sense.

The classical interpretation for the problem is as a partial differential equation in $\Omega \setminus \sulut{D_0 \cap D_\infty}$ with different boundary conditions on $\doo D_0$ and $\doo D_\infty$, as follows:
\begin{remark}
\label{rmk:pde-bdy}
Let the connected components of $D_\infty$ be written as~$C$.

The conductivity equation can be reformulated as
\begin{equation}
\label{eq:traditional_direct}
\begin{cases}
\dive(\sigma\abs{\nabla u}^{p-2}\nabla u)=0 & \text{ in }\Omega\setminus(\ol D_0\cup\ol D_\infty)\\
u=f & \text{ on }\partial\Omega\\
\sigma\abs{\nabla u}^{p-2}\partial_\nu u=0 & \text{ on }\partial D_0\\
\text{for each component~$C$ of~$D_\infty$}
&
\begin{cases}
u|_C = \text{constant}\\
\int_{\partial C}\sigma\abs{\nabla u}^{p-2}\partial_\nu u=0
.
\end{cases}
\end{cases}
\end{equation}
The constant may be different for different components and the constant values depend on the boundary data~$f$.
\end{remark}

The inverse problem of studying inclusions with finite or zero conductivity is well-known~\cite{Isakov:1988,Colton:Kirsch:1996,Cakoni:Colton:2006,Kirsch:Grinberg:2008,Harrach:2013,Ikehata:1998,Potthast:2005,Ikehata:1999:jan,Nakamura:Uhlmann:Wang:2005,Harrach:Ullrich:2013}.
The region of infinite conductivity has been investigated in a few papers; Gorb and Novikov~\cite{Gorb:Novikov:2012} consider the behaviour of the forward problem when there are two inclusions close to each other for the $p$-conductivity equation for $2 \leq p \in \N$.
We are not aware of any results concerning the inverse problem for~$p \neq 2$, but there are results when~$p = 2$.
Br\"uhl~\cite[section 4.3.1]{Bruhl:1999} and Schmitt~\cite[section 2.2.2]{Schmitt:2010} have used the factorization method to detect perfectly conducting inclusions.
Ramdani and Munnier~\cite{Munnier:Ramdani:2015} detect the infinitely conductive bodies from the Dirichlet to Neumann map in two dimensional domain for the linear conductivity equation.
Their method is based on geometry in the complex plane, in particular Riemann mappings.
Friedman and Vogelius~\cite{Friedman:Vogelius:1989} have shown that one can recover the location and scale of a finite number of small inclusions with zero or infinite conductivity in an inhomogeneous background from the DN map.

\subsection{Solutions of Wolff} \label{sec:wolff}
We call certain special solutions for the $p$-Laplace equation the Wolff solutions, as they were originally used by Thomas Wolff~\cite[section 3]{Wolff:2007} when investigating the boundary behaviour of $p$-harmonic functions.
Wolff only defined the solutions for $2 < p < \infty$, though his proof works when $p > 3/2$.
Lewis~\cite{Lewis:1988} extended Wolff's work to the case $1<p<2$ by duality arguments.
Salo and Zhong~\cite[section 3]{Salo:Zhong:2012} provided a unified treatment for all $1 < p < \infty$.
\begin{lemma}
Let $\rho, \rho^{\perp} \in \R^d$ satisfy $\abs{\rho} = \abs{\rho^{\perp}} = 1$ and $\rho \cdot \rho^{\perp} = 0$. Define $h \colon \R^d \to \R$ by $h(x) = e^{-\rho \cdot x}w(\rho^\perp \cdot x)$, where the function~$w$ solves the differential equation
\begin{equation}\label{eq:wolff}
w''(s) + V(w,w')w = 0
\end{equation}
with
\begin{equation}
V(a,b) = \frac{(2p-3)b^2+(p-1)a^2}{(p-1)b^2 + a^2}.
\end{equation}
The function~$h$ is then $p$-harmonic.

Given any initial conditions $(a_0,b_0) \in \R^2 \setminus \{(0,0)\}$ there exists a solution $w \in C^\infty(\R)$ to the differential equation~\eqref{eq:wolff} which is periodic with period $\lambda_p > 0$, satisfies the initial conditions $(w(0),w'(0)) = (a_0,b_0)$, satisfies $\int_0^{\lambda_p} w(s) \der s = 0$, and furthermore there exist constants~$c$ and~$C$ depending on $a_0,b_0,p$ such that for all $s \in \R$ we have
\begin{equation}\label{eq:wolff_new}
C > w(s)^2+w'(s)^2 > c > 0.
\end{equation}
\end{lemma}
When $p = 2$ the equation~\eqref{eq:wolff} is simpler: $w + w'' = 0$.
The solution in this special case is a linear combination of sine and cosine functions.
The solutions for general $p$ behave in a similar way: they oscillate with mean zero and period $\lambda_p > 0$.

\subsection{Dirichlet to Neumann map} \label{sec:dn}

The strong formulation for the Dirichlet to Neumann map, which we simply call the strong Dirichlet to Neumann map, associates Neumann boundary values to Dirichlet data~$f$:
\begin{equation}
\Lambda_\sigma (f) = \sigma \abs{\nabla u}^{p-2}\nabla u \cdot \nu |_{\doo \Omega}
\end{equation}
The strong Dirichlet to Neumann map is defined pointwise and hence requires the conductivity~$\sigma$ to be extendable pointwise to $\doo \Omega$ and also requires $\nabla u$ to be defined pointwise on the boundary.
Due to these difficulties the weak formulation of the Dirichlet to Neumann map, or the weak Dirichlet to Neumann map, is used.
For properties of the Dirichlet to Neumann map with constant conductivity we refer to Hauer~\cite{Hauer:2015}.
For $p$-conductivity equation the weak Dirichlet to Neumann map was introduced by Salo and Zhong~\cite{Salo:Zhong:2012} and used in the included articles~\cite{Brander:Kar:Salo:2015:thesis,Brander:2016:thesis}.
In article~\cite{Brander:Ilmavirta:Kar:2015:thesis} we extend the weak Dirichlet to Neumann map to conductivities that include regions of zero or infinite conductivity, as follows.

We assume~$\Omega$ and~$\sigma$ to be as in theorem~\ref{thm:direct}.
Let $X=W^{1,p}(\Omega)/W^{1,p}_0(\Omega)$ and~$X'$ be  its dual.
The DN map $\Lambda_\sigma\colon X\to X'$ is defined by
\begin{equation}
\ip{\Lambda_\sigma f}{g} = \int_\Omega\sigma\abs{\nabla\bar f}^{p-2}\nabla\bar f\cdot\nabla\bar g \;\der x,
\end{equation}
where $\bar f\in W^{1,p}\sulut{\Omega}$ is any minimizer of the energy functional~$E$ with boundary values $f\in X$ and $\bar g\in W^{1,p}(\Omega)$ is an extension of $g\in X$ with $\nabla \bar g = 0$ in $D_\infty$.
Since $\ol D_\infty\cap\partial\Omega=\emptyset$, there always exists such an extension~$\bar g$.

In article~\cite{Brander:2016:thesis} we give sufficient conditions for recovering the strong Dirichlet to Neumann map from the weak one.
We state the result here with slightly relaxed assumptions:
\begin{lemma}\label{lemma:weaktostrong}
Suppose that $\Omega$ has $C^3$-smooth boundary, boundary values~$f \in C^3(\doo \Omega)$ and that $\nabla \sigma$ is H{\"o}lder-continuous and the conductivity bounded and strictly positive in~$\ol \Omega$.
Then we can recover the pointwise values of the strong Dirichlet to Neumann map
\begin{equation}
\sigma(x_0) \left|\nabla u(x_0)\right|^{p-2} \nabla u (x_0) \cdot \nu(x_0)
\end{equation}
from the weak Dirichlet to Neumann map
\begin{equation}
\left\langle \Lambda_\sigma(f), g \right\rangle = \int_\Omega \sigma |\nabla u|^{p-2} \nabla u \cdot \nabla \bar g \der x.
\end{equation}
\end{lemma}

\section{$p$-Calder\'on's problem}

The inverse problem of Calder\'on was posed in the setting of the $p$-conductivity equation by Salo and Zhong~\cite{Salo:Zhong:2012}.
See section~\ref{subsec:boundary} for their results.

There are some differences between the 2-Calder\'on's problem and $p$-Calder\'on's problem.
In particular, many results in the linear case are based on reducing the conductivity equation to Schrödinger's equation
\begin{equation}
\sulut{-\Delta + q}u = 0
\end{equation}
with potential
\begin{equation}
q = \frac{\Delta \sqrt \sigma}{\sqrt \sigma}.
\end{equation}
It is not clear if a similar reduction is available with the $p$-conductivity equation.

A standard method in investigating Calder\'on-type problems for equations with weak non-linearities~\cite{Sun:1996} is taking the G\^ateaux derivative of the Dirichlet to Neumann map at constant boundary values.
The method does not work in case of the $p$-Calder\'on's problem.
We follow the presentation of Salo and Zhong~\cite[appendix]{Salo:Zhong:2012}.
To see the problem, interpret~$a \in \R$ as Dirichlet boundary data, suppose $t > 0$ and $f \in W^{1,p}\sulut{\Omega}/W^{1,p}_0\sulut{\Omega}$.
Write the solution of $p$-conductivity equation with boundary data~$g$ as $u_g$.
Observe that $u_{a+tf} = a + u_{tf} = a + tu_f$.
Then the G\^ateaux derivative is the limit as $t \to 0$ of
\begin{equation}
\begin{split}
&\frac{1}{t} \sulut{\Lambda_\sigma (a+tf) - \Lambda_\sigma (a)} \\
&= \frac{1}{t}\sigma \abs{\nabla \sulut{a+tu_f}}^{p-2} \nabla \sulut{a+tu_f}\cdot \nu - 0  \\
&= t^{p-2}\Lambda_\sigma (f).
\end{split}
\end{equation}
In particular, the G\^ateaux derivative does not even exist when $p < 2$.
This calculation does not provide any new information, unlike in the article of Sun~\cite{Sun:1996}.

Unique continuation results often play a role in Calder\'on's problem.
The results are much more restricted in the case of the $p$-conductivity equation~\cite{Granlund:Marola:2014,Guo:Kar:2015}.
We are not aware of any work on the related Runge approximation property for the~$p$-Laplace equation.

\section{Boundary determination}\label{subsec:boundary}

In an article published in 2012, Salo and Zhong~\cite{Salo:Zhong:2012} showed how to recover conductivity on the boundary of a domain with reasonable regularity assumptions on the conductivity and the boundary of the domain.
Their proof is similar to that of Brown~\cite{Brown:2001}, and Brown and Salo~\cite{Brown:Salo:2006}, who only considered the linear situation where $p=2$.
\begin{theorem} \label{thm:sz}
Suppose $\Omega \subset \R^d$ is a bounded domain with $C^1$ boundary, $d \geq 2$, and the conductivity~$\sigma$ is continuous at a point~$x_0 \in \doo \Omega$.
Then the weak Dirichlet to Neumann map determines $\sigma(x_0)$.
\end{theorem}
They proved the theorem by first using complex-valued and then real-valued boundary values -- complex geometric optics solutions and Wolff solutions (defined in section~\ref{sec:wolff}), respectively, multiplied by a cutoff function focused at $x_0$.
Since results using only real-valued boundary values are stronger than those that use complex-valued functions, further work has used real-valued boundary values exclusively.

In particular, Salo and Zhong use the following lemma with $g = \sigma$.
They do not explicitly mention the lemma in their paper.
\begin{lemma} \label{lemma:sz}
Suppose that $\Omega \subset \R^d$, $d \geq 2$, is a bounded open set with $C^1$-smooth boundary, and $x_0 \in \doo \Omega$.

Then there exists an explicit sequence of real-valued Dirichlet boundary values~$f_N$ such that the corresponding solutions of the $p$-conductivity equation satisfy
\begin{equation}
\lim_{N \to \infty} q(N)\int_\Omega g(x) \abs{\nabla u_N}^p \der x = g(x_0),
\end{equation}
where $q(N)$ is an explicit scaling constant and $g$ is any function that is continuous in a neighbourhood of the boundary point~$x_0$.
\end{lemma}

We improve the result in~\cite{Brander:2016:thesis} by determining the gradient of conductivity on the boundary:
\begin{theorem} \label{thm:boundary}
Suppose $\Omega \subset \R^d$, $d \geq 2$, is a bounded open set with $C^3$-smooth boundary.
If the  conductivity~$\sigma$ is positive and of class~$C^2\sulut{\ol \Omega}$, then one can recover $\nabla \sigma|_{\doo \Omega}$ from the weak Dirichlet to Neumann map.
\end{theorem}
The extra regularity assumptions on conductivity~$\sigma$ and boundary~$\doo \Omega$ are needed to establish a Rellich-type identity (theorem~\ref{thm:rellich}) and to derive the strong Dirichlet-to-Neumann map from the weak one in lemma~\ref{lemma:weaktostrong}.
The assumptions in the theorem~\ref{thm:boundary} presented above are stronger than in the paper, as are the assumptions of the Rellich-type identity (theorem~\ref{thm:rellich}).

A similar Rellich-type identity was used in the linear situation, where $p = 2$, by Brown, García, and Zhang~\cite[appendix]{Garcia:Zhang:2012}.
\begin{theorem}[Rellich identity, \cite{Brander:2016:thesis}] \label{thm:rellich}
Suppose $\Omega \subset  \R^d$, $d \geq 2$, is a bounded domain with $C^2$-smooth boundary.
Suppose that a function~$u$ solves the $p$-conductivity equation with $1 < p< \infty$ and conductivity $\sigma \in C^{1} \sulut{\ol \Omega}$ bounded above zero.
Let $\alpha \in \R^d$.
Then
\begin{equation}
\begin{split}
\int_\Omega &\sulut{\alpha \cdot \nabla \sigma} \abs{\nabla u}^p \der x \\
&= \int_{\doo \Omega} \alpha \cdot \nu \; \sigma \abs{\nabla u}^p \der S(x) - p \int_{\doo \Omega}  \alpha \cdot \nabla u \; \sigma \abs{\nabla u}^{p-2} \doo_\nu u \; \der S(x).
\end{split}
\end{equation}
\end{theorem}
The identity gives the integral $\int_\Omega g \abs{\nabla u}^p \der x$ in lemma~\ref{lemma:sz} with $g = \alpha \cdot \nabla \sigma$ a formulation that only depends on known quantities: $\sigma$ and $\nabla u$ on the boundary.
The conductivity is known due to results of Salo and Zhong, while the gradient can be recovered from the strong Dirichlet to Neumann map, see lemma~\ref{lemma:weaktostrong}.

\section{Detecting inclusions}

An inclusion is a set where conductivity is significantly higher or lower than in the rest of the domain.
We only consider inclusions embedded in constant background conductivity, which we take to be one.

That is, let $D \subset \Omega$ be an open set so that $\sigma (x) = 1$ whenever $x \notin D$ and let exactly one of the following be true for all $x \in D$:
\begin{enumerate}
\item $\sigma(x) = 0$
\item $0 < c < \sigma(x) < C < 1$
\item $1 < c < \sigma(x) < C < \infty$
\item $\sigma(x) = \infty$.
\end{enumerate}
Here $c$ and $C$ are constants that do not depend on the spatial variable~$x$.

In case of finite non-zero conductivity we show that for $1<p<\infty$ we can detect the convex hull of the inclusion; this is done in article~\cite{Brander:Kar:Salo:2015:thesis}.
In case of zero or infinite conductivity  and $1<p<\infty$ we can detect some convex set $K$, which is a superset of the convex hull of the inclusion; this is done in article~\cite{Brander:Ilmavirta:Kar:2015:thesis}.
We use the enclosure method of Ikehata~\cite{Ikehata:1999:oct:enc}.
We next formulate the results as theorems:
\begin{theorem}
\label{thm:enc_inf}
Suppose that $\Omega \subset \R^d$ is open and bounded with a priori known constant conductivity~$\sigma$ outside an obstacle $D = D_0 \cup D_\infty$.
Suppose $\sigma \colon \Omega \to [0,\infty]$ is measurable.
Suppose either~$D_0$ or~$D_\infty$ is empty, the sets~$\doo \Omega$, $\ol{D_0}$ and~$\ol{D_\infty}$ are pairwise disjoint, and the inclusion~$D$ has Lipschitz boundary.

Then we can, from knowledge of the Dirichlet to Neumann map, find a set~$D'$, which is a superset of the convex hull of~$D$.
Furthermore, we can detect whether~$D_0$ or~$D_\infty$ is non-empty.
\end{theorem}
The proof of the previous theorem can be found in the article~\cite{Brander:Ilmavirta:Kar:2015:thesis}.
By the boundary reconstruction result of Salo and Zhong, theorem~\ref{thm:sz}, if the domain~$\Omega$ is of class~$C^1$, then we don't need to know the value of the background conductivity a priori.
\begin{theorem}
Suppose $\Omega\subset\mathbb{R}^d, d\geq 2$, is a bounded domain and the inclusion $D\subset\Omega$ is a bounded open set with Lipschitz boundary.
Furthermore assume that the measurable conductivity $\sigma$ has a jump discontinuity along the interface~$\doo D$; that is,
$\sigma(x) := 1 + \sigma_D(x)\chi_D(x)$, where $\sigma_D$ is bounded away from zero and either positive and bounded from above or negative and bounded from below by a constant greater than minus one, and where $\chi_D$ is the characteristic function of $D$.

Then the convex hull of $D$ can be recovered from the Dirichlet to Neumann map.
Further, we know if $\sigma_D$ is positive or negative.
\end{theorem}
The previous theorem is proven in article~\cite{Brander:Kar:Salo:2015:thesis}.

The enclosure method relies on an indicator function, which is a difference of Dirichlet to Neumann maps that take particlar Wolff solutions as input.

\begin{definition}
The indicator function~$I$ is defined as 
\begin{equation}
I\sulut{t,\rho,\rho^\perp,\tau} = \tau^{d-p}\ip{\sulut{\Lambda_\sigma - \Lambda_1}f}{f},
\end{equation}
where  $\tau > 0$,
\begin{equation}
f\sulut{x,t,\rho,\rho^\perp,\tau} = e^{\tau\sulut{\rho \cdot x - t}}w\sulut{\tau \rho^\perp \cdot x}
\end{equation} 
is a Wolff solution (considered in the space of boundary values)
and $\sulut{\Lambda_\sigma - \Lambda_1}f = \Lambda_\sigma(f) - \Lambda_1(f)$.
\end{definition}

Note that $\abs{\nabla u(x)} \to 0$ as $\tau \to \infty$ when $\rho \cdot x < t$ and $\abs{\nabla u(x)}$ blows up as $\tau \to \infty$ when $\rho \cdot x > t$.
Here $u$ is the Wolff solution with boundary values~$f$.
Hence, at least heuristically, the parameters $t$ and $\rho$ divide the space so that one half-space has plenty of energy while the other has very little.
The idea is that the indicator function reveals if the inclusion intersects the high-energy half space or is disjoint from it.
Combining this information for all parameters $t$ and $\rho$ yields the convex hull of the inclusion.

When the inclusion is of finite non-zero conductivity (as investigated in article~\cite{Brander:Kar:Salo:2015:thesis}) we have the following theorem:
\begin{theorem} \label{thm:enclosure_estimate}
There exist $c,C > 0$ such that
\begin{align}
c < \abs{I\sulut{\sup_{x \in D} \sulut{x \cdot \rho},\rho, \rho^\perp, \tau}} < C \tau^d
\end{align}
for $\tau \gg 1$.
\end{theorem}
We only need to consider the case $t =\sup_{x \in D} x \cdot \rho$, since the identity
\begin{equation}
  I\sulut{t,\rho, \rho^\perp, \tau} = e^{2\tau(t_0-t)} I\sulut{t_0,\rho, \rho^\perp, \tau}
\end{equation}
with $t_0 = \sup_{x \in D} \sulut{x \cdot \rho}$ shows that the behaviour of the indicator function is exponentially increasing or decreasing as a function of~$\tau$ for other values of~$t$.
The previous theorem is proven with the help of a monotonicity inequality:
\begin{theorem}\label{thm:mono}
If $ 0 < c < \sigma_0, \sigma_1 \in L^{\infty}\sulut{\Omega}$ and $1 < p < \infty$, and if $f \in W^{1,p}(\Omega)$, then 
\begin{align}
(p-1) & \int_{\Omega} \frac{\sigma_0}{\sigma_1^{1/(p-1)}} \sulut{\sigma_1^{\frac{1}{p-1}} - \sigma_0^{\frac{1}{p-1}}} \abs{\nabla u_0}^p \der x \\
& \leq ((\Lambda_{\sigma_1} - \Lambda_{\sigma_0})f, f) 
 \leq \int_{\Omega} (\sigma_1 - \sigma_0) \abs{\nabla u_0}^p \der x,
\end{align}
where $u_0 \in W^{1,p}(\Omega)$ solves $\dive(\sigma_0 \abs{\nabla u_0}^{p-2} \nabla u_0) = 0$ in $\Omega$ with $u_0|_{\partial \Omega} = f$.
\end{theorem}
Monotonicity inequalities have been used in the linear case~\cite{Harrach:Ullrich:2013}.

In case of inclusion with zero or perfect conductivity, as discussed in article~\cite{Brander:Ilmavirta:Kar:2015:thesis}, we only detect when the inclusion intersects the half-space of high energy, hence possibly detecting too large a set.
This means that we only have the lower bound in theorem~\ref{thm:enclosure_estimate}.

\section{Other interior results}
Very recently, Guo, Kar and Salo~\cite{Guo:Kar:Salo:2016} used the monotonicity inequality, theorem~\ref{thm:mono}, to show injectivity for the Dirichlet to Neumann map under a monotonicity assumption.
Their result show injectivity in two dimensions for Lipschitz conductivities.
In higer dimensions they need an additional assumption;  one of the conductivities must be almost constant.

\newpage

\bibliographystyle{plain}
\bibliography{math}


\cleardoublepage
\thispagestyle{empty}

~\vspace{15em}

\begin{center}

[I]
\\[2em]
{\bf Enclosure method for the $p$-Laplace equation}, \\
T.~Brander, M.~Kar, M.~Salo, \\
Inverse Problems  31(4):045001 (2015).
\\[3em]
\textcopyright\ IOP Publishing Ltd.\\
Published article: \url{http://dx.doi.org/10.1088/0266-5611/31/4/045001} \\
Preprint: \url{http://arxiv.org/abs/1410.4048}

\end{center}

\cleardoublepage



\cleardoublepage
\thispagestyle{empty}

~\vspace{15em}

\begin{center}

[II]
\\[2em]
{\bf Calder\'on problem for the $p$-Laplacian: First order derivative of conductivity on the boundary}, \\
T.~Brander, \\
Proceedings of American Mathematical Society 144 (2016), 177-189. 
\\[3em]
\textcopyright\ American Mathematical Society.\\
Published article: \url{http://dx.doi.org/10.1090/proc/12681} \\
Preprint: \url{http://arxiv.org/abs/1403.0428}

\end{center}

\cleardoublepage



\cleardoublepage
\thispagestyle{empty}

~\vspace{15em}

\begin{center}

[III]
\\[2em]
{\bf Superconductive and insulating inclusions for non-linear conductivity equations}, \\
T.~Brander, J.~Ilmavirta, M.~Kar, \\
preprint.

\end{center}

\cleardoublepage

\includepdf[pages={-}]{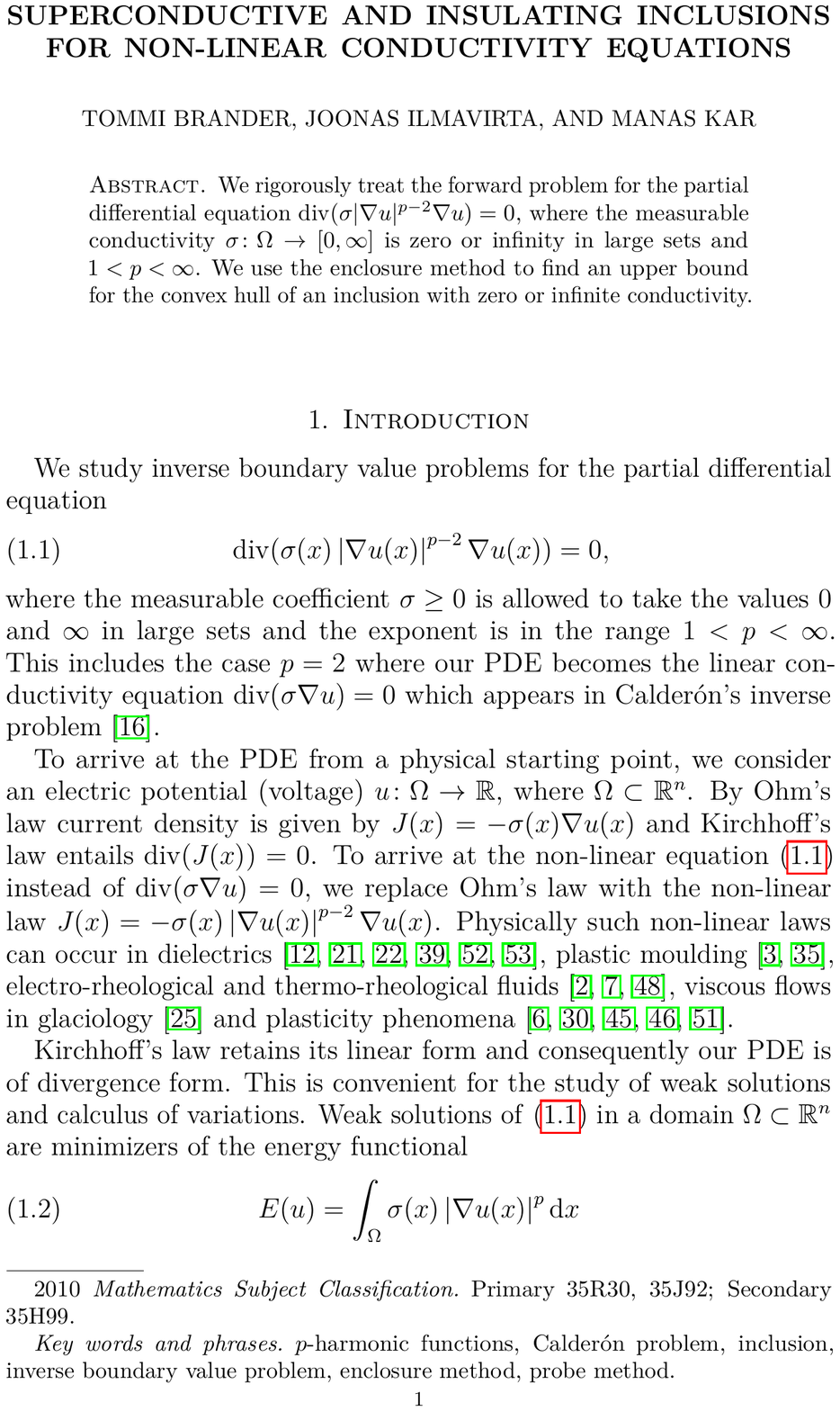}




\end{document}